\documentclass[12pt, intlimits]{amsart}
\usepackage{amsmath,amsthm,amscd,amssymb}
\usepackage{latexsym}

\allowdisplaybreaks
\newcommand{\Hmm}[1]{\leavevmode{\marginpar{\tiny%
$\hbox to 0mm{\hspace*{-0.5mm}$\leftarrow$\hss}%
\vcenter{\vrule depth 0.1mm height 0.1mm width \the\marginparwidth}%
\hbox to 0mm{\hss$\rightarrow$\hspace*{-0.5mm}}$\\\relax\raggedright #1}}}

\def\bea{\begin{eqnarray}}
\def\eea{\end{eqnarray}}

\def\R{\mathbb R}
\def\C{\mathbb C}

\def\half{\frac12}
\def\eps{\epsilon}

\newcommand{\nopar}{\vspace{2mm}\noindent}

\newcommand{\fint}{{\int}{\kern-1em \diagup}}
\newcommand{\pref}[1]{{\rm(\ref{#1})}}

\newcommand{\bbR}{{\mathbb{R}}}

\newcommand{\calH}{{\mathcal H}}
\newcommand{\calS}{{\mathcal S}}
\newcommand{\calF}{{\mathcal F}}

\newcommand{\f}{\frac}

\newcommand{\rank}{\text{\rm{rank}}}

\newcommand{\beq}{\begin{equation}}
\newcommand{\eeq}{\end{equation}}
\newcommand{\ba}{\begin{align}}
\newcommand{\ea}{\end{align}}

\newcommand{\lang}{\left\langle}
\newcommand{\rang}{\right\rangle}

\newtheorem{theorem}{Theorem}
\newtheorem*{theorem*}{Theorem}

\newtheorem{lemma}[theorem]{Lemma}
\newtheorem*{lemma*}{Lemma}
\newtheorem{cor}[theorem]{Corollary}
\newtheorem*{corollary*}{Corollary}
\theoremstyle{definition}

\newtheorem*{definition*}{Definition}

\theoremstyle{remark}
\newtheorem{remark}{Remark}

\def\la{\langle}
\def\ra{\rangle}
\def\nn{\nonumber}
\def\beeq{\begin{equation}}
\def\eneq{\end{equation}}

\begin{document}

\title[Dispersive estimates]{Dispersive estimates for Schr\"{o}dinger operators in the presence of a resonance and/or
an eigenvalue at zero energy in dimension three: I}

\author{M. Burak Erdo\smash{\u{g}}an and Wilhelm Schlag}
\thanks{This work was done in June of 2004,
while the first author visited Caltech and he  wishes to thank that institution for its hospitality and support. 
The first author was partially supported by the NSF grant DMS-0303413.
The second author was partially supported by a Sloan fellowship and the NSF grant DMS--0300081.}

\address{Department of Mathematics \\
University of Illinois \\
Urbana, IL 61801, U.S.A.}

\address{Department of Mathematics \\
Caltech \\
Pasadena, CA 91125, U.S.A.}

\date{\today}

\maketitle

\section{Introduction}

\nopar
Consider the Schr\"{o}dinger operator  $H=-\Delta+V$ in $\R^3$, where $V$ is a real-valued potential.
Let $P_{ac}$ be the
orthogonal projection onto the absolutely continuous subspace of $L^2(\R^3)$ which is determined by $H$.
In \cite{JSS}, \cite{Y1}, \cite{RodSch}, \cite{goldbergschlag} and \cite{goldberg}, $L^1(\R^3)\rightarrow L^\infty(\R^3)$ dispersive estimates
for the time evolution $e^{itH}P_{ac}$ were investigated under various decay assumptions on the potential $V$ and the
assumption that zero is neither an eigenvalue nor a resonance of~$H$. Recall that zero energy is a resonance iff
there is $f\in L^{2,-\sigma}(\R^3)\setminus L^2(\R^3)$ for all  $\sigma>\half$ so that $Hf=0$. Here $L^{2,-\sigma}=\langle x\rangle^{\sigma}L^2$
are the usual weighted $L^2$ spaces and $\langle x\rangle:=(1+|x|^2)^{\frac12}$.

\nopar
In the first part of this paper we investigate dispersive estimates when there is a resonance at energy zero.
It is well-known, see  Rauch~\cite{Rauch}, Jensen and Kato~\cite{JenKat}, and Murata~\cite{Mur}, that the decay in
that case is~$t^{-\half}$. Moreover, these authors derived expansions of the evolution into inverse powers of time
in weighted $L^2(\R^3)$ spaces. Here, we obtain such expansions with respect to the $L^1\to L^\infty$ norm, albeit only in terms of the
powers $t^{-\half}$ and $t^{-\frac32}$.
Our results
will require decay of the form
\bea\label{decay}
|V(x)|\leq C \langle x\rangle^{-\beta},
\eea
for some $\beta>0$. Our goal was brevity rather than optimality. In particular, it was not our intention to
obtain the minimal value of $\beta$, and our results can surely
be improved in that regard.
Our first result is for the case when zero is only a resonance, but not an eigenvalue.

\begin{theorem}\label{T:scalar1}
Assume that $V$ satisfies \pref{decay} with $\beta>12$ and assume that there is a resonance at energy zero
but that zero is not an eigenvalue. Then there is a time dependent rank one operator $F_t$ (see  \pref{f_t}  below) 
such that
$$
\left\|e^{itH} P_{ac}-t^{-1/2} F_t  \right\|_{1\rightarrow\infty}\leq C  t^{-3/2},
$$
for all $t>0$ and  $F_t$ satisfies
\bea\label{fbound}
\sup_t\left\|F_t\right\|_{L^1\rightarrow L^\infty}<\infty, \qquad \limsup_{t\to\infty}\left\|F_t\right\|_{L^1\rightarrow L^\infty} >0.
\eea
\end{theorem}

\nopar
The following case allows for any combination of resonances and/or eigenvalue at energy zero. It is important to note that
the $t^{-\frac32}$ bound is destroyed by an eigenvalue at zero, even if zero is not a resonance and even after projecting the
zero eigenfunction away (this was discovered by Jensen and Kato~\cite{JenKat}).

\begin{theorem}
\label{T:scalar2}
Assume that $V$ satisfies \pref{decay} with $\beta>12$ and assume that there is a resonance at energy zero
and/or that zero is an eigenvalue. Then there is a time dependent operator $F_t$ 
such that
$$
\sup_t\left\|F_t\right\|_{L^1\rightarrow L^\infty}<\infty, \qquad \left\|e^{itH} P_{ac}-t^{-1/2}F_t  \right\|_{1\rightarrow\infty}\leq C  t^{-3/2}.
$$
\end{theorem}

\nopar In all cases, the operators $F_t$ can be given explicitly,
and they can of course be extracted from our proofs. The methods
of this paper also apply to matrix Schr\"{o}dinger operators, as
considered for example in Cuccagna~\cite{Cuc} or~\cite{Sch}.
Details of this will be given elsewhere.

\section{Scalar case}
\nopar
Let
$K_{\lambda_0}$ be the operator with kernel
$$
K_{\lambda_0}(x,y)=\f{1}{\pi i}\int_0^\infty e^{it\lambda^2}\lambda \chi_{\lambda_0}(\lambda)  [R_V^+(\lambda^2)-R_V^-(\lambda^2)](x,y)
d\lambda,
$$
where
$$
R_V^\pm(\lambda^2)= R_V(\lambda^2\pm i0)=(H-(\lambda^2 \pm i0))^{-1}
$$
is the perturbed resolvent. By the limiting absorption principle, these boundary values are bounded operators on
weighted $L^2$-spaces, see e.g.~\cite{agmon}.
Here $\chi$ is an even smooth function supported in $[-1,1]$ and $\chi(x)=1$ for $|x|<1/2$;
$\chi_{\lambda_0}(\lambda)=\chi(\lambda/\lambda_0)$.
The high energies were studied in \cite{goldbergschlag}:
\begin{theorem}\cite{goldbergschlag} Assume that $V$ satisfies \pref{decay} with some $\beta>3$, then for any $\lambda_0>0$
$$
\left\|e^{itH} P_{ac}-K_{\lambda_0}  \right\|_{1\rightarrow\infty}\leq C_{\lambda_0} t^{-3/2}.
$$
\end{theorem}

\nopar
Hence, in the proof of Theorem~\ref{T:scalar1} and Theorem~\ref{T:scalar2}, it suffices to consider the
operator $K_{\lambda_0}$ for some $\lambda_0$.
One can rewrite the kernel of $K_{\lambda_0}$ as
\bea \label{kxy}
K_{\lambda_0}(x,y)=\f{1}{\pi i}\int_{-\infty}^\infty e^{it\lambda^2}\lambda \chi_{\lambda_0}(\lambda)   R_V((\lambda+i0)^2)(x,y)
d\lambda,
\eea
Note that $R((\lambda+i0)^2)(x,y)$ is not an even function of $\lambda$; rather, we have
$$
R_V((\lambda+i0)^2)(x,y)=\overline{R_V\left((-\lambda+i0)^2\right)(x,y)}.
$$

\subsection{Resolvent expansions at zero energy}

\nopar
In this section, following \cite{jensennenciu}, we obtain resolvent expansions at the threshold
$\lambda=0$ in the presence of a resonance.  This is of course similar to Jensen and Kato~\cite{JenKat}, but
we prefer to work with the $L^2$-based approach from~\cite{jensennenciu}.
For $j=0,1,2,...,$ let $G_j$ be the operator with  the kernel
$$
G_j(x,y)=\f{1}{4\pi j!}|x-y|^{j-1}.
$$
Recall that for each $J=0,1,2,...$,
\bea\label{resexp}
R_0(\lambda^2)=\sum_{j=0}^{J} (i\lambda)^j G_j+o(\lambda^J), \text{ as }\lambda\rightarrow 0.
\eea
This  expansion is valid
in the space, $HS_{L^{2,\sigma}\rightarrow L^{2,-\sigma}}$, of Hilbert-Schmidt operators between $L^{2,\sigma}$ and $L^{2,-\sigma}$
for $\sigma>\max((2J+1)/2,3/2)$.

\nopar
Let $U(x)=1$ if $V(x)\geq 0$ and $U(x)=-1$ if $V(x)<0$, $v=|V|^{1/2}$ and $w=vU$. We have
$$
V=Uv^2=wv.
$$
We use the symmetric resolvent identity, valid for $\Im\lambda>0$:
\bea\label{res_exp}
R_V(\lambda^2)=  R_0(\lambda^2)-R_0(\lambda^2)vA(\lambda)^{-1}vR_0(\lambda^2),
\eea
where
\begin{align}\label{alamla}
A (\lambda ) &= U + v R_0 (\lambda^2) v=(U + v G_0 v)+\lambda \f{v[R_0 (\lambda^2) - G_0 ]v}{\lambda}\\
&=:A_0+\lambda A_1 (\lambda).\nonumber
\end{align}
$A_1(\lambda)$ has the kernel
\begin{align*}
A_1(\lambda)(x,y)&=\f{1}{\lambda}v(x)\f{e^{i\lambda|x-y|}-1}{4\pi|x-y|}v(y),\\
|A_1(\lambda)(x,y)|&\leq\f{1}{4\pi}|v(x)|\;|v(y)|.
\end{align*}
Therefore, $A_1(\lambda)\in HS:=HS_{L^2\rightarrow L^2}$ provided $\langle x \rangle^{\f{3}{2}+} v(x) \in L^\infty$.
Also note that
$$
A_1 (0) = iv G_1 v = \f{i\alpha}{4\pi} P_v,\;\;\alpha=\|V\|_1,
$$
where $P_v$ is the orthogonal projection onto  span$(v)$. It is important to realize
that $A(\lambda)$ has a natural meaning for $\lambda\in\R$ via the limit $\R+i0$.

\nopar
First, we consider the expansions of   $ A(\lambda)^{-1}$ for $\lambda$ close to
zero as in \cite{jensennenciu}. The following  lemma (Corollary 2.2 in \cite{jensennenciu})
is our main tool. Note the similarity between \pref{az-1} and the symmetric resolvent identity.

\begin{lemma}\cite{jensennenciu}\label{L:jen-nen}
Let $F\subset \C\setminus\{0\}$ have zero as an accumulation point. Let $A(z)$, $z\in F$,
be a family of bounded operators of the form
$$A(z)=A_0+z A_1(z)$$
with $A_1(z)$ uniformly bounded as $z\rightarrow 0$. Suppose that $0$ is an isolated point of the spectrum of
$A_0$, and let $S$ be the corresponding Riesz projection. Assume that $\rank(S)<\infty$.
Then for sufficiently small $z\in F$ the operators
$$
B(z):=\frac{1}{z}(S-S(A(z)+S)^{-1} S)
$$
are well-defined and bounded on $\calH$. Moreover, if $A_0=A_0^*$, then they are
 uniformly bounded as $z\rightarrow 0$. The operator $A(z)$ has a bounded inverse in $\calH$ if and only if
$B(z)$ has a bounded inverse in $S\calH$, and in this case
\bea\label{az-1}
A(z)^{-1}=(A(z)+S)^{-1}+\frac{1}{z}(A(z)+S)^{-1}SB(z)^{-1}S(A(z)+S)^{-1}.
\eea
\end{lemma}
\begin{proof}
It is a standard fact that
\[ {\mathrm Ran}(S)\supset \bigcup_{n=1}^\infty \ker(A_0^n).\]
By our assumption $\rank(S)<\infty$ we have equality here, and $(A_0+S)^{-1}$ has a bounded inverse.
Hence, $A(z)+S$ also has a bounded inverse for small~$z$, as can be seen from the usual Neuman series.
Therefore, $B(z)$ is well-defined for small $z$ and bounded. Moreover, if $A_0$ is self-adjoint, then
\[ S-S(A_0+S)^{-1}S=0\]
which implies that $B(z) =O(1)$ as $z\to0$.
Suppose $B(z)$ is invertible on $S\calH$. Denote the right-hand side of \eqref{az-1} by $T(z)$. Then
\begin{align*}
 T(z)A(z)&=A(z)T(z)\\
 &=  I + \frac{1}{z}SB(z)^{-1}S(A(z)+S)^{-1}-S(A(z)+S)^{-1} \\
&\quad -\frac{1}{z}S(A(z)+S)^{-1}SB(z)^{-1}S(A(z)+S)^{-1} \\
&= I + \frac{1}{z}SB(z)^{-1}S(A(z)+S)^{-1}-S(A(z)+S)^{-1} \\
&\quad-\frac{1}{z}(S-zB(z))B(z)^{-1}S(A(z)+S)^{-1}=I.
\end{align*}
Conversely, suppose that $A(z)$ is invertible. Define
\[ D(z):= z(S+SA(z)^{-1}S)=z(A(z)+S)[A(z)^{-1}-(A(z)+S)^{-1}](A(z)+S). \]
Then
\begin{align*}
B(z)D(z)&=D(z)B(z) \\
&= S+SA(z)^{-1}S-S(A(z)+S)^{-1}S-S(A(z)+S)^{-1}SA(z)^{-1}S \\
&= S+S(A(z)+S)^{-1}SA(z)^{-1}S-S(A(z)+S)^{-1}SA(z)^{-1}S=S,
\end{align*}
so that $D(z)$ is the inverse of $B(z)$ on $S\calH$.
\end{proof}

\nopar
Note that $A_0$ as in \eqref{alamla} is a compact perturbation of $U$ and that the essential spectrum of $U$ is contained in $\{-1,1\}$.
Moreover, $A_0$ is self adjoint.
Therefore, $0$ is an isolated point of the spectrum of
$A_0$ and $\dim(\ker_{A_0})<\infty$. Let $S_1$ be the corresponding Riesz projection. Since $A_0$ is self adjoint, $S_1$ is the orthogonal projection
onto the kernel of $A_0$ and we have
\bea\label{s1commute}
S_1=(A_0+S_1)^{-1}S_1=S_1(A_0+S_1)^{-1}.
\eea
\begin{remark} \label{absbound} By the resolvent identity we have
$$(A_0+S_1)^{-1}=U-(A_0+S_1)^{-1}(vG_0v+S_1)U.$$
Since $|V(x)|\lesssim \langle x\rangle^{-3-}$ and $S_1$ is a finite rank operator, we have $(vG_0v+S_1)U\in HS$,
and hence
$(A_0+S_1)^{-1}$ is the sum of $U$ and
a Hilbert-Schmidt operator. Therefore, the operator with kernel
$|(A_0+S_1)^{-1}(x,y)|$ is bounded in $L^2$.
This remark will be useful below when we consider dispersive estimates.
\end{remark}

\nopar
We choose $\lambda_0>0$ sufficiently small so that $A (\lambda) + S_1$ is invertible for   $|\lambda|<\lambda_0$.
Using Lemma~\ref{L:jen-nen}, we see that, for $|\lambda|<\lambda_0$,
$A(\lambda)$ is  invertible if and only if
$$
m(\lambda) = \f{S_1 - S_1 \left(A (\lambda) + S_1\right)^{-1}S_1}{\lambda}
$$
is invertible on $S_1L^2$ and in this case
\bea
\label{aters}
A(\lambda)^{-1}  =  (A(\lambda) + S_1 )^{-1} +\f{1}{\lambda}
 (A(\lambda) + S_1 )^{-1} S_1 m(\lambda)^{-1}
S_1  (A(\lambda)+S_1 )^{-1}.
\eea
If $\lambda_0$ is sufficiently small, then
\bea\nonumber
\left(A (\lambda) + S_1\right)^{-1}
 =(A_0 + S_1)^{-1} + \sum_{k=1}^\infty (-1)^k \lambda^k (A_0 + S_1)^{-1}
\left(A_1 (\lambda) (A_0+S_1)^{-1}\right)^k.
\eea
Plugging this into the definition of $m(\lambda)$ and using \pref{s1commute},
we obtain
\begin{align*}
m(\lambda) &= S_1A_1(\lambda)S_1 + \sum^\infty_{k=1} (-1)^k \lambda^k S_1
\left(A_1 (\lambda)(A_0+S_1)^{-1}\right)^{k+1} S_1\\
&=m(0)+\lambda m_1(\lambda),
\end{align*}
where
\begin{align}
m (0) &= S_1A_1(0)S_1=
  \f{i\alpha}{4\pi} S_1 P_v S_1, \nn\\
m_1 (\lambda)&= S_1
\f{A_1 (\lambda) - A_1 (0)}{\lambda}S_1
+\sum_{k=1}^{\infty} (-1)^k \lambda^{k } S_1
\left(A_1(\lambda) (A_0 + S_1)^{-1}\right)^{k+1}S_1.\label{eq:m1def}
\end{align}
If $m(0)$ is invertible in $S_1L^2$, then we can invert $m(\lambda)$ for small $\lambda$ using Neuman series
 and hence obtain an expansion for $A(\lambda)^{-1} $. Since $m(0)$ has rank one, this can only occur if $\rank(S_1)=1$.

\nopar
Otherwise, let $S_2:S_1L^2\to S_1L^2$ be the orthogonal projection onto the kernel of $m(0)$ where the latter operates on~$S_1L^2$.
As above, $m(\lambda)+S_2$ is invertible in $S_1L^2$ for
$|\lambda|<\lambda_0$ (we choose a smaller $\lambda_0$ if necessary), and
\bea\label{s2commute}
S_2=S_2(m(0)+S_2)^{-1}=(m(0)+S_2)^{-1}S_2.
\eea
Lemma~\ref{L:jen-nen} asserts that
$m(\lambda)$ is invertible on $S_1L^2$ if and only if
$$
b(\lambda)= \f{S_2 - S_2 \left(m (\lambda) + S_2\right)^{-1} S_2}{\lambda}
$$
is invertible on $S_2L^2$ and
\bea\label{mters}
 m(\lambda)^{-1} = \left(m (\lambda)+S_2\right)^{-1}
+ \f{1}{\lambda} \left(m (\lambda ) + S_2 \right)^{-1}
S_2 b (\lambda)^{-1} S_2
\left(m (\lambda) + S_2\right)^{-1}.
\eea
Note that
\bea\nonumber
\left(m(\lambda) + S_2\right)^{-1}
 =(m(0) + S_2)^{-1} + \sum_{k=1}^\infty (-1)^k \lambda^k (m(0) + S_2)^{-1}
\left(m_1(\lambda) (m(0)+S_2)^{-1}\right)^k.
\eea
Plugging this into the definition of $b(\lambda)$ and using \pref{s2commute},
we obtain
\begin{align*}
b(\lambda) &= S_2 m_1(\lambda)S_2 + \sum_{k=1}^\infty (-1)^k\lambda^k S_2\left(m_1(\lambda) (m(0)+S_2)^{-1}\right)^{k+1}S_2\\
&=:b(0)+\lambda b_1(\lambda),
\end{align*}
where $b(0)=S_2m_1(0)S_2$ and
\beeq
\label{eq:b1def}
b_1(\lambda)  = \f{S_2 [m_1(\lambda)-m_1(0)]S_2}{\lambda} +
\f{1}{\lambda} \sum_{k=1}^\infty (-1)^k\lambda^k S_2\left(m_1(\lambda) (m(0)+S_2)^{-1}\right)^{k+1}S_2.
\eneq
A simple calculation using \pref{resexp} (with $J=2$) and $S_2S_1=S_1S_2=S_2$
shows that
\begin{equation}
\label{eq:b(0)}
b(0) =-S_2 v G_2 v S_2.
\end{equation}
Since $G_2\in HS_{L^{2, \sigma}\rightarrow L^{2,-\sigma}}$   for $\sigma>5/2$,
we have $b(0)\in HS$  if $|V(x)|\lesssim   \langle x \rangle^{-5-\varepsilon}$.

\nopar
Below, we characterize the spaces $S_1L^2$, $S_2L^2$ and also prove that $b(0)$ is always invertible on $S_2L^2$.
Therefore, for small $\lambda$, $b (\lambda)$ is invertible.
This proves that $A (\lambda)$ is invertible for $0<|\lambda|<\lambda_0$. Using \pref{aters} and \pref{mters}, we obtain
\begin{align}\label{expansion}
  A (\lambda)^{-1}  &  = \Gamma_1(\lambda) \\
&+ \f{1}{\lambda} \Gamma_1(\lambda) S_1
\Gamma_2(\lambda) S_1 \Gamma_1(\lambda) \nonumber \\
& + \f{1}{\lambda^2} \Gamma_1(\lambda) S_1 \Gamma_2(\lambda)
S_2 b(\lambda)^{-1} S_2 \Gamma_2(\lambda) S_1 \Gamma_1(\lambda),\nonumber
\end{align}
where
$$
\Gamma_1(\lambda)=\left(A (\lambda) + S_1\right)^{-1},\text{ and }
 \Gamma_2(\lambda)=\left(m (\lambda) + S_2 \right)^{-1}.
$$
Note that this formula is also valid in the case $S_2=0$.
\begin{lemma}\label{L:S1} Assume $|V(x)| \lesssim \langle x \rangle^{-3-\varepsilon}$. Then
$f \in S_1L^2\backslash\{0\}$ if and only if $f=wg$ for some $g \in L^{2,-\f{1}{2}-}\backslash\{0\}$ such that
\bea \label{Hg=0}
-\Delta g + Vg =0 \ \text{in}\ \calS^\prime.
\eea
\end{lemma}
\begin{proof}
First recall that, for $g\in L^{2,-\f{1}{2}-}$,  (\ref{Hg=0}) holds  if and only if
$$
(I+G_0V)g=0,
$$
see Lemma 2.4 in \cite{JenKat}.
Suppose $f \in S_1L^2\backslash\{0\}$.  Then
\begin{align}
\nonumber
(A_0f)(x)&=U(x) f(x) + \f{v(x)}{4\pi} \int \f{v(y)f(y)}{|x-y|}   \;dy  =0  \\
\Rightarrow  \;\; & f(x) + \f{w(x)}{4\pi}
\int \f{v(y)f(y)}{|x-y|} \;dy  =0.
\label{fsifir}
\end{align}
Let
\bea\label{g_def}
g(x)=-\f{1}{4\pi}\int_{\R^3} \f{v(y)f(y)}{|x-y|} \;dy.
\eea
Note that $g\in L^{{2,-\f{1}{2}-}}$ and $f(x)=w(x)g(x)$ for each $x$. Moreover,
(\ref{Hg=0}) holds since
$$
g(x)=-\f{1}{4\pi}\int_{\R^3} \f{v(y)f(y)}{|x-y|} \;dy=-\f{1}{4\pi}\int_{\R^3} \f{V(y)g(y)}{|x-y|} \;dy=-[G_0V g](x).
$$

\nopar
Conversely, assume $f=wg$ for some $g$ as in the hypothesis. Then $f \in L^{2,1+}$ and
\begin{align*}
A_0f(x)&=U(x) f(x) + \f{v(x)}{4\pi} \int \f{v(y)f(y)}{|x-y|}   \;dy   \\
& = v(x)g(x)+\f{v(x)}{4\pi} \int \f{V(y)g(y)}{|x-y|}   \;dy  =v(I+G_0V)g=0.
\end{align*}
Note that since $g$ is not identically zero, $Vg\neq 0$, and hence $f\neq 0$.
\end{proof}

\nopar
By Lemma~\ref{L:S1}, we see that   $f\in S_1L^2$ implies   $f\in L^{2, 1+}$.
\begin{lemma}\label{L:S2} Assume $|V(x)| \lesssim \langle x \rangle^{-3-\varepsilon}$. Then
$f \in S_2L^2\backslash\{0\}$ if and only if $f=wg$ for some $g \in L^{2}\backslash\{0\}$ such that
$$-\Delta g + Vg =0 \ \text{in}\ \calS^\prime.$$
\end{lemma}
\begin{proof} Suppose $f \in S_2L^2\backslash\{0\}$. Note that $S_2L^2\subset S_1L^2$ and,
 by Lemma~\ref{L:S1}, we have $f=wg$ for some $g\in L^{2,-\f{1}{2}-}\backslash\{0\}$ such that
$-\Delta g + Vg =0 \ \text{in}\ \calS^\prime$.  By the definition of $S_2$, we have
$$
S_1P_v f=0.
$$
Note that
$S_1P_v f =0$ if and only if
$$
S_1 v  = 0  \text{ or }
P_v f = 0.
$$
In the first case, $S_2=S_1$ and $P_v f=0$ for any $f\in S_2L^2$.  We have the same conclusion
in the second case. Thus,
$$
\int_{\R^3} v(y)f(y)\;dy=0.
$$
Using this and (\ref{g_def}), we obtain
$$
g(x)  = -\f{1}{4\pi} \int_{\bbR^3}
\left[\f{1}{|x-y|} - \f{1}{1 + |x|}\right]
v  (y) f  (y) \;dy \in L^{2,\f{1}{2} - }.
$$
This is because
\bea\label{inek}
\left|\f{1}{|x-y|} - \f{1}{1 + |x|}\right| \leq
\f{1 + |y|}{|x-y|(1 + |x|)}
\eea
and $f \in L^{2,1+}$.

\nopar
Conversely, assume $f=wg$ for some $g$ as in the hypothesis. Then
$$
g(x)  = -\f{1}{4\pi} \int_{\bbR^3}
\left[\f{1}{|x-y|} - \f{1}{1 + |x|}\right]
v  (y) f  (y) \;dy - \f{1}{4\pi(1+|x|)} \int_{\bbR^3}
v  (y) f  (y) \;dy.
$$
By \pref{inek} the first summand is in $L^2$. Therefore
$$
\left[\int_{\R^3} v(y)f(y)\;dy\right] \f{1}{1+|\cdot|}\in L^2(\R^3).
$$
Thus, $\int v(y)f(y)\;dy=0$ and $f\in S_2L^2\backslash \{0\}$.
\end{proof}
\begin{lemma}\label{L:b0}
Assume $|V(x)| \lesssim \langle x \rangle^{-5-\varepsilon}$. Then,
as an operator in $S_2L^2$, the kernel of $b(0)$ is trivial.
\end{lemma}
\begin{proof}
Assume that for some $f\in S_2L^2$, $b(0)f=0$, i.e.,
$$
\langle G_2 v f, vf\rangle=0.
$$
From the proof of Lemma \pref{L:S2}, we have
$$
\int_{\R^3} f(y)v(y) dy =0.
$$
Using this and \pref{resexp} (with $J=2$), we obtain
\begin{align*}
0&= \langle G_2 vf, vf\rangle\\
&=\lim_{\lambda \to 0}
\lang\f{R_0 (\lambda^2)-G_0}{\lambda^2} v f , v f \rang\\
&=\lim_{\lambda\to 0} \f{1}{\lambda^2} \int\left((\xi^2 + \lambda^2)^{-1} - \xi^{-2}\right)
\widehat{v f}(\xi) \overline{\widehat{v f}}(\xi)d\xi\\
&=\lim_{\lambda\to 0} \int \f{1}{\xi^2(\xi^2+\lambda^2)}
 |\widehat{v f}(\xi)|^2 d\xi  \\
&= \int \f{|\widehat{vf}(\xi)|^2}{\xi^4} d \xi\; \text{ (by the Monot.\ Conv.\ Thm.)}\\
&= \la R_0(0)vf,R_0(0)vf\ra  \qquad \Rightarrow \widehat{vf}= 0 \Rightarrow vf = 0.
\end{align*}
Using this in \pref{fsifir}, we obtain $f=0$.
\end{proof}

\subsection{Dispersive estimate when zero is not an eigenvalue}
\label{subsec:S2=0}

\nopar
In this section, we prove Theorem~\ref{T:scalar1}. When zero is not an eigenvalue, $S_2=0$ and \pref{expansion} reduces to
\bea\label{S2=0}
A(\lambda)^{-1}  = (A(\lambda) + S_1)^{-1}
 + \f{1}{\lambda} (A(\lambda) + S_1)^{-1} S_1
 m(\lambda)^{-1} S_1 (A(\lambda) + S_1)^{-1},
\eea
where
\begin{align}
\left(A (\lambda) + S_1\right)^{-1}
&=(A_0 + S_1)^{-1} + \sum_{k=1}^\infty (-1)^k \lambda^k (A_0 + S_1)^{-1}
\left[A_1 (\lambda) (A_0+S_1)^{-1}\right]^k\nn\\
&=:(A_0+S_1)^{-1}+\lambda E_1(\lambda),\nn \\
m(\lambda)^{-1}&=m(0)^{-1}+\sum_{k=1}^\infty (-1)^k \lambda^k  m(0)^{-1}
\left[m_1(\lambda) m(0)^{-1}\right]^k \label{eq:m0inv}\\
&=:m(0)^{-1}+\lambda E_2(\lambda).\nn
\end{align}
Thus, using \pref{s1commute}, we obtain
\begin{align}
\label{S2=0expan}
A (\lambda)^{-1}  & =  \frac{1}{\lambda}S_1  m(0)^{-1} S_1  \\
& + \left(A (\lambda) + S_1\right)^{-1} \nonumber
\\
 &+   E_1(\lambda) S_1
 m(\lambda)^{-1} S_1 \left(A (\lambda) + S_1\right)^{-1}\nonumber \\
 &+   \left(A (\lambda) + S_1\right)^{-1} S_1
 E_2(\lambda) S_1 \left(A (\lambda) + S_1\right)^{-1} \nonumber\\
 &+   \left(A (\lambda) + S_1\right)^{-1} S_1
 m(\lambda)^{-1} S_1 E_1(\lambda)  \nonumber \\
& =: \f{1}{\lambda}S+E(\lambda).\nonumber
\end{align}
Note that $S$ is a rank one operator. Plugging (\ref{S2=0expan}) into (\ref{res_exp}), we have
\begin{align*}
R_V (\lambda^2) &=-\f{1}{\lambda} R_0(\lambda^2)vSv R_0(\lambda^2)\\
&+R_0 (\lambda^2) - R_0(\lambda^2)vE(\lambda)vR_0(\lambda^2).
\end{align*}
Using this in (\ref{kxy}), we get
$$
K_{\lambda_0}(x,y)=  K_1(x,y) +  K_2(x,y) - K_3(x,y),
$$
where
\begin{align}
K_1(x,y)&  =\f{-i}{16\pi^3 }\int_{-\infty}^\infty \int_{\R^6} e^{it\lambda^2} \chi_{\lambda_0}(\lambda)
\f{e^{i\lambda(|x-u_1|+|y-u_2|)}}{|x-u_1||y-u_2|}v(u_1)S (u_1,u_2)v(u_2)  d u_1 d u_2   d\lambda,\nn \\
K_2(x,y)& = \int_{-\infty}^\infty e^{it\lambda^2}\lambda \chi_{\lambda_0}(\lambda) R_0 (\lambda^2)(x,y)  d \lambda \nn \\
K_3(x,y)& = \int_{-\infty}^\infty e^{it\lambda^2}\lambda \chi_{\lambda_0}(\lambda)[R_0(\lambda^2)vE(\lambda)vR_0(\lambda^2)](x,y)  d \lambda.
\label{eq:K3def}
\end{align}
First,  we deal with $K_1$. Note that
\begin{align}\label{k1}
&K_1(x,y)= \\
  \f{-i}{16\pi^3}&\int_{\R^6}\int_{-\infty}^\infty  e^{it\lambda^2} \chi_{\lambda_0}(\lambda)\f{\cos(\lambda(|x-u_1|+|y-u_2|))}{|x-u_1||y-u_2|}v(u_1)S (u_1,u_2)v(u_2)  d u_1 d u_2   d\lambda. \nonumber
\end{align}
We have
\begin{align} \label{fourier}
 \int_{-\infty}^\infty  t^{1/2}e^{it\lambda^2} \chi_{\lambda_0}(\lambda) \cos(\lambda a)   d\lambda &=
\int_{-\infty}^\infty \left(t^{1/2}e^{it(\cdot)^2}\right)^\vee(u)
{\left(\chi_{\lambda_0}(\cdot) \cos(\cdot a) \right)}^{\wedge}(u)du\\
&=c\int_{-\infty}^\infty e^{iu^2/4t}(\widehat{\chi_{\lambda_0}}(u+a)+\widehat{\chi_{\lambda_0}}(u-a))du \nonumber\\
&=c\int_{-\infty}^\infty e^{i(u^2+a^2)/4t} \cos(\f{ua}{2t})\widehat{\chi_{\lambda_0}}(u)du\nonumber\\
&=c\int_{-\infty}^\infty e^{i(u^2+a^2)/4t}  \widehat{\chi_{\lambda_0}}(u)du\nonumber\\
&+c\int_{-\infty}^\infty e^{i(u^2+a^2)/4t} (\cos(\f{ua}{2t})-1)\widehat{\chi_{\lambda_0}}(u)du\nonumber\\
&=:C_1(t,a)+C_2(t,a).\nonumber
\end{align}
Using this in (\ref{k1}), we obtain
\begin{align*}
K_1(x,y)&=\f{-it^{-1/2}}{16\pi^3}\int_{\R^6} \f{C_1(t,|x-u_1|+|y-u_2|)}{|x-u_1||y-u_2|}v(u_1)S (u_1,u_2)v(u_2)  d u_1 d u_2   \\
&+\f{-it^{-1/2}}{16\pi^3}\int_{\R^6} \f{C_2(t,|x-u_1|+|y-u_2|)}{|x-u_1||y-u_2|}v(u_1)S (u_1,u_2)v(u_2)  d u_1 d u_2 \\
&=:K_{11}(x,y)+K_{12}(x,y).
\end{align*}
Note that
$$
|C_2(t,a)|\leq c \f{|a|}{t}.
$$
Thus,
\begin{align}\label{kbirki}
|K_{12}(x,y)|&\leq c t^{-3/2} \int_{\R^6} \left(\f{1}{|x-u_1|}+\f{1}{|y-u_2|}\right) |v(u_1)||v(u_2)||S (u_1,u_2)|   d u_1 d u_2 \\
&\lesssim t^{-3/2} \left[\left\|\f{v(\cdot)}{|x-\cdot|}\right\|_2 +\left\|\f{v(\cdot)}{|y-\cdot|}\right\|_2\right]\| |S| \|_{2\rightarrow 2} \|v\|_2
 \nonumber\\
 &\lesssim t^{-3/2}. \nonumber
\end{align}
The last inequality follows from the fact that  $S$ is a rank one operator and the following
calculation which holds for $v\in L^2\cap L^\infty$;
\begin{align}
\label{vcalc}
\left\| \f{|v(\cdot)| }{|x-\cdot| }\right\|_2^2& =\int_{|x-u|<1} \f{|v(u)|^2 }{|x-u|^2 } du + \int_{|x-u|>1} \f{|v(u)|^2 }{|x-u|^2 } du
\\
&\lesssim \int_{|u|<1} \f{1}{|u|^2 } du + \int_{\R^3}  |v(u)|^2   du \lesssim 1.
\nonumber
\end{align}
Now, we consider $K_{11}$. Note that
$$
C_1(t,a)=e^{ia^2/4t}h(t),
$$
where $h(t)$ is a smooth function which converges to $c$ as $t$ tends to $\infty$. We have
\begin{align}
K_{11}(x,y)&=\f{-ih(t)}{16\pi^3t^{1/2}}\int_{\R^6}
\f{e^{i|x-u_1|^2/4t}e^{i|y-u_2|^2/4t}}{|x-u_1||y-u_2|}v(u_1)S (u_1,u_2)v(u_2)  d u_1 d u_2 \nonumber  \\
& - \f{ih(t)}{16\pi^3t^{1/2}}\int_{\R^6}
\f{e^{\f{i(|x-u_1|+|y-u_2|)^2}{4t}}-e^{\f{i(|x-u_1|^2+|y-u_2|^2)}{4t}}}{|x-u_1||y-u_2|}v(u_1)
S(u_1,u_2)v(u_2)  d u_1 d u_2 \nonumber \\
&=:t^{-1/2}F_t(x,y)+K_{112}(x,y). \label{f_t}
\end{align}
Since $S$ is a rank one operator, for each $t$,   $F_t$ is a rank one operator. Also note that
by a calculation similar to \pref{kbirki}, we obtain $\sup_{t,x,y}|F_t(x,y)|\lesssim 1$. Finally, $F_t\ne0$
for all $t$, and $\lim_{t\to\infty}F_t$ exists in the weak sense and does not vanish:
\[
\lim_{t\to\infty} \langle F_t f,g\rangle = \frac{-i c}{16\pi^3} \int_{\R^{12}} \frac{f(x)\bar{g}(y)}{|x-u_1||y-u_2|} v(u_1)S(u_1,u_2)v(u_2)\, du_1du_2\,dxdy
\]
for any $f,g\in\calS$.
By a similar calculation, the term $K_{112}$
is dispersive since
$$
\left| e^{i(|x-u_1|+|y-u_2|)^2/4t}-e^{i(|x-u_1|^2+|y-u_2|^2)/4t}\right|\lesssim \f{|x-u_1||y-u_2|}{t}.
$$
$K_2$ is the low energy part of the free evolution and hence it is dispersive.
The rest of this section is devoted to the proof of
\begin{equation}
\label{eq:K3est}
\sup_{x,y}|K_3(x,y)|\lesssim t^{-3/2}.
\end{equation}
Denote
$$\f{d}{d\lambda}\left(\chi_{\lambda_0}(\lambda)R_0(\lambda^2)vE(\lambda)vR_0(\lambda^2)\right)$$
by $\calF_{x,y}(\lambda)$.
By integration by parts we obtain
$$
K_3(x,y)= \f{1}{2it}\int_{-\infty}^\infty e^{it\lambda^2}\calF_{x,y}(\lambda) d\lambda.
$$
Using Parseval's formula, we obtain
\beeq
\label{eq:K3begin}
K_3(x,y)=\f{c}{t^{3/2}}\int_{-\infty}^\infty e^{i\xi^2/4t} \widehat{\calF_{x,y}}(\xi)d\xi.
\eneq
Thus, it suffices to prove that
\bea\label{fl1}
\sup_{x,y} \| \widehat{\calF_{x,y}}\|_{L^1}<\infty.
\eea
Recall that
$$
\calF_{x,y}(\lambda)=\int_{\R^6}\f{d}{d\lambda}\left[\chi_{\lambda_0}(\lambda)E (\lambda)(u_1,u_2)v(u_1)v(u_2)
\f{e^{i\lambda(|x-u_1|+|y-u_2|)}}{|x-u_1||y-u_2|}\right]du_1du_2.
$$
Let us  concentrate on the term where the derivative hits $\chi_{\lambda_0}(\lambda)E (\lambda)$ (the term
where the derivative hits the exponential is similar):
$$
\tilde{\calF}_{x,y}(\lambda)=\int_{\R^6}[\chi_{\lambda_0}(\lambda)E (\lambda)]^\prime(u_1,u_2)v(u_1)v(u_2)
\f{e^{i\lambda(|x-u_1|+|y-u_2|)}}{|x-u_1||y-u_2|}du_1du_2.
$$
Note that
\begin{align*}
\|\widehat{\tilde{\calF}_{x,y}}(\xi)\|_{L^1}&=\int_{-\infty}^\infty \left| \int_{\R^6}
\widehat{(\chi_{\lambda_0}E )^\prime}(\xi+|x-u_1|+|y-u_2|)(u_1,u_2)
\f{v(u_1)v(u_2)}{|x-u_1||y-u_2|}du_1du_2 \right|d\xi\\
&\leq \int_{\R^6} \int_{-\infty}^\infty \left|\widehat{(\chi_{\lambda_0}E )^\prime}(\xi+|x-u_1|+|y-u_2|)(u_1,u_2)\right|
\f{|v(u_1)||v(u_2)|}{|x-u_1||y-u_2|}  d\xi du_1du_2\\
&=\int_{\R^6} \int_{-\infty}^\infty \left|\widehat{(\chi_{\lambda_0}E )^\prime}(\xi)(u_1,u_2)\right|
\f{|v(u_1)||v(u_2)|}{|x-u_1||y-u_2|}  d\xi du_1du_2\\
&\leq \left\| \f{|v(\cdot)| }{|x-\cdot| }\right\|_2
\left\| \f{|v(\cdot)| }{|y-\cdot| }\right\|_2
\int_{-\infty}^\infty \left\|\left|\widehat{(\chi_{\lambda_0}E )^\prime}(\xi) \right| \right\|_{L^2\rightarrow L^2}d\xi\\
&\lesssim \int_{-\infty}^\infty \left\|\,\left|\widehat{(\chi_{\lambda_0}E )^\prime}(\xi) \right|\, \right\|_{L^2\rightarrow L^2}d\xi.
\end{align*}
The second line follows from Minkowski's inequality and Fubini's theorem, the third line follows from a change of variable,
and the last line follows from the calculation \pref{vcalc}.
Therefore, for $\tilde{\calF}_{x,y}$, (\ref{fl1}) follows from
\beeq
\int_{-\infty}^\infty
\left\|\,\left|\widehat{(\chi_{\lambda_0}E )^\prime}(\xi) \right|\, \right\|_{L^2\rightarrow L^2}d\xi
< \infty.
\label{eq:K3end}
\eneq
We shall use the following elementary lemma.

\begin{lemma}\label{L:convolve}
For each $\lambda\in \R$, let $F_1(\lambda)$ and $F_2(\lambda)$ be bounded operators from $L^2(\R^3)$ to $L^2(\R^3)$
with kernels $K_1(\lambda)$ and $K_2(\lambda)$. Suppose that $K_1, K_2$ both have compact support in $\lambda$ and
that $K_j(\cdot)(x,y)\in L^1(\R)$ for a.e.~$x,y\in\R^3$, as well as
\begin{align}
 &\sup_{\lambda\in\R}\int_{\R^3} |K_1(\lambda)(x_1,x_2)K_2(\lambda)(x_2,x_3)|\, dx_2 <\infty \label{eq:ass1}\\
 &\int_{\R^3}\int_{\R^2} |\widehat{K}_1(\xi)(x_1,x_2)| |\widehat{K}_2(\eta)(x_2,x_3)|\,d\xi d\eta dx_2 < \infty \label{eq:ass2}
\end{align}
for a.e.~$x_1,x_3\in\R^3$.
Let $F(\lambda)=F_1(\lambda)\circ F_2(\lambda)$ with kernel $K(\lambda)$. Then
$$
\int_{-\infty}^\infty \left\|\,\left|\widehat{K}(\xi)\right|\,\right\|_{ 2\rightarrow 2} d\xi
\leq \left[\int_{-\infty}^\infty \left\|\, \left|\widehat{K_1 }(\xi)\right|\,\right\|_{ 2\rightarrow 2}d\xi\right]
\left[\int_{-\infty}^\infty \left\|\, \left|\widehat{K_2 }(\xi)\right|\,\right\|_{ 2\rightarrow 2}d\xi\right].
$$
\end{lemma}
\begin{proof}
By definition, for a.e.~$x_1,x_3\in\R^3$,
\[ K(\lambda)(x_1,x_3) = \int_{\R^3} K_1(\lambda)(x_1,x_2)K_2(\lambda)(x_2,x_3)\, dx_2\]
 and $K(\cdot)(x_1,x_3)\in L^\infty(\R)\cap L^1(\R)$ by \eqref{eq:ass1} and the compact support assumption in~$\lambda$.
Moreover, for a.e.~$\xi\in\R$,
\begin{equation} \hat{K}(\xi)(x_1,x_3) = \int_{\R^3} \int_{-\infty}^\infty \hat{K_1}(\xi-\eta)(x_1,x_2)
\hat{K_2}(\eta)(x_2,x_3) \, d\eta dx_2.
\label{eq:ftF}
\end{equation}
To see this final indentity, denote the right-hand side by $F(\xi;x_1,x_3)$. Then $F(\cdot;x_1,x_3)\in L^1(\R)$ for
a.e.~choice of $x_1,x_3$ by~\eqref{eq:ass2}, and
\begin{align*}
 \int_{-\infty}^\infty e^{2\pi i\xi} F(\xi;x_1,x_3)\, d\xi &= \int_{-\infty}^\infty\int_{\R^3} \int_{-\infty}^\infty e^{2\pi i(\xi-\eta)}\hat{K_1}(\xi-\eta)(x_1,x_2) e^{2\pi i \eta}\hat{K_2}(\eta)(x_2,x_3) \, d\eta dx_2 d\xi \\
& = \int_{\R^3} K_1(\lambda)(x_1,x_2)K_2(\lambda)(x_2,x_3)\, dx_2 d\lambda.
\end{align*}
The final equality sign here follows by Fubini and since $\hat{K_j}(\cdot)(x,y)\in L^1(\R)$ for a.e.~choice of $x,y$ by~\eqref{eq:ass2}.
Hence, \eqref{eq:ftF} holds by uniqueness of the Fourier transform.
The lemma now follows by putting absolute values inside of~\eqref{eq:ftF} and duality.
\end{proof}

\nopar
Note that
$\f{d}{d\lambda}[\chi_{\lambda_0}(\lambda)E (\lambda)]$ is a sum of operators each of which is a composition of  operators
from the list below  (here $\chi(\lambda)$ is a suitably chosen smooth cutoff supported in $[-\lambda_0,\lambda_0]$):
\begin{align*}
F_1(\lambda) &=  \chi(\lambda) (A (\lambda)+S_1)^{-1}, \\
F_2(\lambda) &=  \chi(\lambda) E_1 (\lambda),\\
F_3(\lambda) &=  \chi(\lambda) S_1 m (\lambda)^{-1}S_1, \\
F_4(\lambda) &=  \chi(\lambda) S_1 E_2 (\lambda) S_1,
\end{align*}
and their $\lambda$ derivatives. Moreover, we leave it to the reader to check that for each of the combinations that
contribute to $E(\lambda)$ the hypotheses of Lemma~\ref{L:convolve} are fulfilled.
Therefore, in light of Lemma~\ref{L:convolve}, the following lemma completes the analysis of $K_3$.
\begin{lemma}\label{L:hammal}
For each of the operators $F_j$, $j=1,2,3,4$ above,
$$
\int_{-\infty}^\infty \left\| \left|\widehat{F_j }(\xi)\right|\right\|_{ 2\rightarrow 2}d\xi < \infty.
$$
The same statement is valid for their $\lambda$ derivatives, too.
\end{lemma}
\begin{proof}
We omit the analysis of $F_1$ and $F_3$. Recall that
\begin{align*}
F_2(\lambda)&=\chi(\lambda)E_1 (\lambda)=\chi(\lambda)  \f{\left(A  (\lambda) + S_1\right)^{-1}
-(A_0 + S_1)^{-1}}{\lambda} \\
&=  \chi(\lambda) \sum_{k=1}^\infty (-1)^k \lambda^{k-1} (A_0 + S_1)^{-1}
\left[A_1  (\lambda) (A_0+S_1)^{-1}\right]^k.
\end{align*}
Let $\chi_1$ be a smooth cut off function which is equal to $1$ in $[-1,1]$. Note that  the support of
$\chi$ is contained in $[-1,1]$. We have
$$
F_2(\lambda)=  \sum_{k=1}^\infty (-1)^k   \chi(\lambda)\lambda^{k-1} (A_0 + S_1)^{-1}
\left[  \chi_1(\lambda)A_1 (\lambda) (A_0+S_1)^{-1}\right]^k.
$$
Using Lemma~\ref{L:convolve} and Young's inequality, we obtain
\begin{align} \label{F_2}
& \int_{-\infty}^\infty \left\| \left|\widehat{F_2}(\xi)\right|\right\|_{ 2\rightarrow 2}d\xi \leq \\
& \sum_{k=1}^\infty \|\widehat{(\chi(\lambda)\lambda^{k-1})}\|_{L^1} \| |(A_0+S_1)^{-1}| \|_{2\rightarrow 2}^{k+1}
\left[\int_{-\infty}^\infty \| | \widehat{(\chi_1 A_1 )}(\xi)| \|_{2\rightarrow 2} d\xi\right]^k. \nonumber
\end{align}
By Remark~\ref{absbound}, $|(A_0+S_1)^{-1}|$ is bounded on $L^2$. Also note that
\begin{align}
\|\widehat{(\chi(\lambda)\lambda^{k-1})}\|_{L^1}& \lesssim \|(1+|\xi|)\widehat{(\chi(\lambda)\lambda^{k-1})}(\xi)\|_{L^2}\nonumber\\
&\lesssim \|\chi(\lambda)\lambda^{k-1}\|_2 +\| \f{d}{d\lambda}(\chi(\lambda)\lambda^{k-1}) \|_2 \nonumber\\
&\lesssim \lambda_0^k. \label{chilambda}
\end{align}
Below, we prove that
\bea\label{chia1}
\int_{-\infty}^\infty \| | \widehat{(\chi_1 A_1 )}(\xi)| \|_{2\rightarrow 2} d\xi\lesssim 1.
\eea
If $\lambda_0$ is chosen sufficiently small, using  \pref{chilambda} and \pref{chia1} in \pref{F_2} completes the proof of the lemma for $F_2$.
Recall that
\begin{align*}
A_1  (\lambda) (x,y) &=
v(x) \f{e^{i\lambda |x-y|}-1}{4\pi \lambda |x-y|} v(y)\\
&= \f{1}{4\pi i}v(x)v(y) \int_0^1
e^{i\lambda|x-y|b} \;db.
\end{align*}
Therefore,
\begin{align*}
\widehat{(\chi_1 A_1 )}(\xi)(x,y) &=
\f{1}{4\pi i}v(x)v(y)
\int_0^1
\widehat{\chi}(\xi-|x-y|b) \; db.
\end{align*}
Hence by Schur's test, we have
\begin{align} \label{bitartik}
\int_{-\infty}^\infty \left\| \left|  \widehat{\chi A_1 }
(\xi)  \right|\right\|_{2\rightarrow 2} \; d\xi
&\leq \int_{-\infty}^\infty \int_0^1 \sup_x\int_{\R^3}   |v(x)|  | \widehat{\chi}(\xi-|x-y|b)|
|v(y)|  dy \,db \,d\xi\\
&\leq \int_{-\infty}^\infty   \int_0^1 \sup_x \int_{\R^3} |v(x)|
\left(1+|\xi - |x-y|b|\right)^{-1-\varepsilon} |v(y)|dy \,db \,d\xi. \nonumber
\end{align}
Using the inequality
$$
(1+|\alpha-\beta|)\geq (1+|\alpha|)/(1+|\beta|),
$$
we obtain
\begin{align*}
\pref{bitartik} & \leq \int_{-\infty}^\infty     \sup_x \int_{\R^3} |v(x)|
 (1+|\xi|)^{-1-\varepsilon} (1+ |x-y|)^{1+\varepsilon}   |v(y)| dy   \,d\xi \\
 &\lesssim  \left[\sup_x|v(x)|(1+|x|)^{1+\varepsilon}\right]
 \left[\int_{-\infty}^\infty (1+|\xi|)^{-1-\varepsilon}\,d\xi \right]
 \left[ \int_{\R^3}  |v(y)| (1+ | y|)^{1+\varepsilon}  \,  dy\right]\\
 &\leq C_\varepsilon.
\end{align*}
The last line holds for some $\varepsilon>0$ provided $|V(x)|\lesssim \langle x \rangle^{-8-}$.

\nopar
Next, we consider  $F_4$:
$$
F_4(\lambda) = \chi(\lambda) S_1 E_2 (\lambda) S_1
=\chi(\lambda) S_1 \sum_{k=1}^\infty (-1)^k \lambda^{k-1}  m (0)^{-1}
\left[m_1 (\lambda) m (0)^{-1}\right]^k S_1.
$$
Arguing as in the case of $F_2$, it suffices to prove that
\beeq
\label{eq:m1FT}
\int_{-\infty}^\infty \| | \widehat{(\chi_1 m_1 )}(\xi)| \|_{2\rightarrow 2} d\xi\lesssim 1,
\eneq
where $\chi_1$ is a smooth cut off function which is equal to $1$ in the support of $\chi$
(i.e. in $[-\lambda_0,\lambda_0]$)
and which is supported in $[-\lambda_1,\lambda_1]$.
Recall that
$$
m_1(\lambda)=S_1 \f{A_1(\lambda)- A_1 (0)}{\lambda} S_1
+\sum_{j=1}^\infty S_1 (-1)^j \lambda^{j-1}
\left(A_1 (\lambda) (A_0 + S_1)^{-1}\right)^{j+1} S_1.
$$
The second summand can be analysed as above
(here $\lambda_1$ is chosen sufficiently small to guarantee the convergence of the series, and than we choose $\lambda_0$ even smaller).
Now, we consider the first summand. Note that
\begin{align}
A_2 (\lambda)(x,y) &:= \f{A_1(\lambda) - A_1(0)}{\lambda}(x,y) \label{eq:A2def}\\
&= v(x)\f{e^{i\lambda|x-y|}-i\lambda |x-y| -1}{\lambda^2 |x-y|}v(y) \nn\\
&=-v(x)|x-y|v(y)
\int_0^1(1-b)e^{i\lambda |x-y|b} \,db. \nn
\end{align}
Therefore, for any $\varepsilon>0$, we have
\begin{align*}
 \left|\widehat{\chi_1 S_1A_2S_1}(\xi)(x,y)\right|&=
  \left| v(x) |x-y|
v(y)  \int_0^1 (1-b) \widehat{\chi_1}(\xi-b|x-y|) \,db   \right|  \\
&\lesssim |v(x)| |x-y|
|v(y)| \int_0^1 \left(1 + |\xi-b|x-y||)^{-1-\varepsilon}\right)^{-1} \, db\\
&\lesssim |v(x)| |x-y|
|v(y)|  \left(\f{1+|x-y|}{1+|\xi|}\right)^{1+\varepsilon}.
\end{align*}
Using Schur's test, we have
\begin{align*}
\int_{-\infty}^\infty \left\| \left|\widehat{\chi_1 S_1A_2S_1}(\xi) \right|\right\|_{2\rightarrow 2} \,d\xi
&\lesssim  \int_{-\infty}^\infty \f{d\xi}{(1+|\xi|)^{1+\varepsilon}}   \sup_x \int_{\R^3}   |v(x)| |x-y|
  |v(y)|  (1+|x-y|)^{1+\varepsilon} \,dy  \\
&\lesssim    \sup_x \left(|v(x)| (1+|x|^{2+\varepsilon}\right) \int_{\R^3} |v(y)| (1+|y|^{2+\varepsilon})\,dy\\
&\lesssim 1,
\end{align*}
for sufficiently small $\varepsilon$ provided
$|V(x)|  \lesssim \langle x \rangle ^{-10-}$.

\nopar
Next, we  deal with   $\f{d}{d\lambda}F_j(\lambda)$. Once again we omit the analysis of $F_1$ and $F_3$.
Note that

\begin{align*}
\f{d}{d\lambda}F_2(\lambda) = &   \sum_{k=1}^\infty (-1)^k \f{d}{d\lambda}\left(\chi(\lambda) \lambda^{k-1}\right) (A_0 + S_1)^{-1}
\left[A_1 (\lambda) (A_0+S_1)^{-1}\right]^k \\
 +& \sum_{k=1}^\infty (-1)^k \chi(\lambda) \lambda^{k-1} (A_0 + S_1)^{-1}\times\\
&\times\sum_{j=1}^k  [A_1 (\lambda) (A_0+S_1)^{-1} ]^{j-1}  [\f{d}{d\lambda}A_1(\lambda) (A_0+S_1)^{-1}]
 [A_1 (\lambda) (A_0+S_1)^{-1}]^{k-j}
\end{align*}
Arguing as above, it suffices to prove that
\bea\label{chia1prime}
\int_{-\infty}^\infty \| | \widehat{(\chi_1(A_1)^\prime)}(\xi)| \|_{2\rightarrow 2} d\xi\lesssim 1.
\eea
Note that
\begin{align*}
\f{d}{d\lambda}A_1(\lambda) (x,y)
&=   -v(x) \f{e^{i\lambda |x-y|} - i\lambda |x-y| e^{i\lambda|x-y|}-1}
{\lambda^2|x-y|} \; v(y)\\
&= -v(x) \f{e^{i\lambda|x-y|}-i\lambda |x-y|-1}{\lambda^2|x-y|} \;v(y)
+ iv(x) \f{e^{i\lambda|x-y|}-1}{\lambda} v(y)\\
&=-A_2 (\lambda ) + i \widetilde A_1 (\lambda)
\end{align*}
These are similar to the terms treated above. Therefore \pref{chia1prime} holds
provided $|V(x)|  \lesssim \langle x \rangle ^{-10-}$.

\nopar
Finally, we analyze $\f{d}{d\lambda}F_4(\lambda)$. In view of the preceding, it suffices to prove that
\bea\label{a2prime}
\int_{-\infty}^\infty \| | \widehat{(\chi_1(A_2)^\prime)}(\xi)| \|_{2\rightarrow 2} d\xi\lesssim 1.
\eea
We have
\begin{align*}
\f{d}{d\lambda}A_2 (\lambda)(x,y) =& v(x) i \f{e^{i\lambda|x-y|} -1}{\lambda^2} v(y) -
2v (x) \f{e^{i\lambda|x-y|} -i\lambda |x-y| -1} {\lambda^3|x-y|} v(y)\\
 =& -2v (x) \f{e^{i\lambda |x-y| } + \f{1}{2} \lambda^2 |x-y|^2 - i\lambda |x-y| -1}{\lambda^3 |x-y|} v(y)
\\
& + iv (x) \f{-i\lambda |x-y| + e^{i\lambda|x-y|}-1}{\lambda^2}v(y)
\end{align*}
These are treated as before; \pref{a2prime} holds provided $|V(x)|  \lesssim \langle x \rangle ^{-12-}$.
\end{proof}

\subsection{The general case}

We now turn to the proof of Theorem~\ref{T:scalar2}. In view of~\eqref{res_exp}, \eqref{eq:b(0)}, and~\eqref{expansion},
the coefficient of the $\lambda^{-2}$ power in~\eqref{res_exp} equals
\[ R_0(0)v\Gamma_1(0)S_1\Gamma_2(0)S_2 b(0)^{-1} S_2\Gamma_2(0)S_1\Gamma_1(0) vR_0(0)
= -G_0 v S_2 [S_2vG_2vS_2]^{-1} S_2 v G_0.\]

\begin{lemma}
\label{lem:P0}
The operator $G_0 v S_2 [S_2vG_2vS_2]^{-1} S_2 v G_0$ equals the orthogonal pojection in $L^2(\R^3)$ onto the eigenspace of
$H=-\Delta+V$ at zero energy.
\end{lemma}
\begin{proof}
Let $\{\psi_j\}_{j=1}^J$ be an orthonormal basis in ${\rm Ran}(S_2)$.
By Lemmas~\ref{L:S1} and~\ref{L:S2},
\[ \psi_j + wG_0 v\psi_j =0 \quad \forall\,1\le j\le J\]
and we can write $\psi_j=w\phi_j$
for $1\le j\le J$ where $\phi_2\in L^2$, and 
$$\int V\phi_j\, dx= \int v\psi_j\, dx=0.$$
Moreover, the $\{\phi_j\}_{j=1}^J$ are linearly independent and they satisfy
\[ \phi_j + G_0 V \phi_j=0\]
for all $1\le j\le J$.
Since for any $f\in L^2(\R^3)$, $S_2f=\sum_{j=1}^J \la f,\psi_j\ra \psi_j$, we conclude that
\[ S_2 vG_0 f = \sum_{j=1}^J \la f,G_0 v \psi_j\ra \psi_j = - \sum_{j=1}^J \la f,\phi_j\ra \psi_j.\]
Let $A=\{A_{ij}\}_{i,j=1}^J$ denote the matrix of the Hermitian operator
\[ S_2 vG_2v S_2 = \frac{1}{8\pi} S_2v(x)|x-y|v(y)S_2 \]
relative to the basis $\{\psi_j\}_{j=1}^J$. Since $\int_{\R^3} v\psi_j\,dx=0$, the proof of Lemma~\ref{L:b0} shows that
\begin{align*}
A_{ij} &= \la \psi_i, S_2 vG_2v S_2 \psi_j \ra = \la G_0v\psi_i, G_0v\psi_j\ra \\
&= \la  G_0V\phi_i, G_0V\phi_j\ra = \la \phi_i,\phi_j\ra.
\end{align*}
Let
\[ Q:= G_0 v S_2 [S_2vG_2vS_2]^{-1} S_2 v G_0. \]
Then for any $f\in L^2(\R^3)$,
\begin{align*}
Qf &= - \sum_{j=1}^J G_0 v S_2 [S_2vG_2vS_2]^{-1}\psi_j \la f,\phi_j \ra \\
&= -\sum_{i,j=1}^J G_0 v S_2 \psi_i(A^{-1})_{ij} \la f,\phi_j \ra = \sum_{i,j=1}^J \phi_i(A^{-1})_{ij} \la f,\phi_j \ra.
\end{align*}
In particular,
\[ Q\phi_k = \sum_{i,j=1}^J \phi_i(A^{-1})_{ij} \la \phi_k,\phi_j \ra = \sum_{i,j=1}^J \phi_i(A^{-1})_{ij} A_{jk} = \phi_k\]
for all $1\le k\le J$. The conclusion is that ${\rm Ran}\,Q={\rm span}\{\phi_j\}_{j=1}^J$, and that $Q={\rm Id}$ on
${\rm Ran}\,Q$. Since $Q$ is Hermitian, it is the orthogonal projection onto ${\rm span}\{\phi_j\}_{j=1}^J$, as claimed.
\end{proof}

\noindent
This has the following simple and standard consequence for the spectral measure.
\begin{cor}
\label{cor:P0}
Let $-\infty<\lambda_{N}<\lambda_{N-1}<\ldots<\lambda_1<\lambda_0\le 0$ be the finitely many
eigenvalues of $H=-\Delta+V$.
Let $P_{\lambda_j}$ denote the orthogonal projection in $L^2(\R^3)$ onto the
 eigenspace of $H$ corresponding to the eigenvalue $\lambda_j$.
Then
\begin{equation}
\label{eq:repform}
 e^{itH} = \sum_{j=0}^N e^{it\lambda_j} P_{\lambda_j} + \frac{1}{2\pi i} \int_0^\infty e^{it\lambda} [R_V^+(\lambda)-R_V^{-}(\lambda)]\, d\lambda.
\end{equation}
Moreover,
\begin{equation}
\label{eq:root}
 R_V^+(\lambda)-R_V^{-}(\lambda) =  O_{L^2}(\lambda^{-\half})
\end{equation}
as $\lambda\to 0+$ so that the integral in~\eqref{eq:repform} is absolutely convergent at $\lambda=0$.
\end{cor}
\begin{proof}
Start from the expression
\[
 e^{itH} =  \frac{1}{2\pi i} \int_0^\infty e^{it\lambda} [R_V(\lambda+i\eps)-R_V(\lambda-i\eps)]\, d\lambda,
\]
which is valid for all $\eps>0$ (via the spectral theorem, for example).
The formula~\eqref{eq:repform} follows by passing to the limit $\eps\to0$. Indeed, the projections
arise as Cauchy integrals
\[ P_{\lambda_j}\,\frac{1}{2\pi i}\oint_{\gamma_j} \frac{dz}{z-\lambda_j} = P_{\lambda_j} \]
where $\gamma_j$ is a small circle surrounding $\lambda_j$. We need to invoke Lemma~\ref{lem:P0} in case $\lambda_0=0$, since
it determines the coefficient of the $z^{-1}$ singularity in the asymptotic expansion of the resolvent.
Once we subtract that singularity, what remains is $O(|z|^{-\half})$, as claimed.
\end{proof}

\nopar
The point of Lemma~\ref{lem:P0} and Corollary~\ref{cor:P0} is
really to prove~\eqref{eq:root}, since~\eqref{eq:repform} is of course obvious.
One can also deduce Lemma~\ref{lem:P0} from the proof of the Corollary starting from~\eqref{eq:repform},
since the most singular power $z^{-1}$ must lead to the projection onto the eigenspace. However, we have chosen
to give these direct proofs.

\begin{proof}[Proof of Theorem~\ref{T:scalar2}]
In view of \eqref{expansion},
\begin{align*}
  A (\lambda)^{-1}  &  = \Gamma_1(\lambda) \\
&+ \f{1}{\lambda} \Gamma_1(\lambda) S_1
\Gamma_2(\lambda) S_1 \Gamma_1(\lambda) \nonumber \\
& + \f{1}{\lambda^2} [\Gamma_1(\lambda) S_1 \Gamma_2(\lambda)
S_2 b(\lambda)^{-1} S_2 \Gamma_2(\lambda) S_1 \Gamma_1(\lambda)-S_2b(0)^{-1}S_2]\\
& + \f{1}{\lambda^2} S_2 b(0)^{-1} S_2.
\end{align*}
Inserting this into \eqref{res_exp} leads to
\begin{align}
R_V(\lambda^2) &=  R_0(\lambda^2)-R_0(\lambda^2)v\Gamma_1(\lambda) vR_0(\lambda^2) \nn \\
& - \f{1}{\lambda} R_0(\lambda^2)v \Gamma_1(\lambda) S_1
\Gamma_2(\lambda) S_1 \Gamma_1(\lambda) vR_0(\lambda^2) \label{eq:S22}\\
& - \f{1}{\lambda^2} R_0(\lambda^2)v [\Gamma_1(\lambda) S_1 \Gamma_2(\lambda)
S_2 b(\lambda)^{-1} S_2 \Gamma_2(\lambda) S_1 \Gamma_1(\lambda)-S_2b(0)^{-1}S_2] vR_0(\lambda^2) \label{eq:S23}\\
& - \f{1}{\lambda^2} (R_0(\lambda^2)-G_0)v S_2 b(0)^{-1} S_2 vR_0(\lambda^2) - \f{1}{\lambda^2} G_0 v S_2 b(0)^{-1} S_2 v(R_0(\lambda^2)-G_0)
\label{eq:S24}\\
& -\f{1}{\lambda^2} P_0. \label{eq:S25}
\end{align}
The three terms up to and including \eqref{eq:S22} have already been covered in Subsection~\ref{subsec:S2=0}.
Indeed, the only difference here is that we need to incorporate $S_2$ into the expression~\eqref{eq:m0inv}:
\begin{align*}
m(\lambda)^{-1}&=(m(0)+S_2)^{-1}+\sum_{k=1}^\infty (-1)^k \lambda^k  (m(0)+S_2)^{-1}
\left[m_1(\lambda) (m(0)+S_2)^{-1}\right]^k \\
&=:(m(0)+S_2)^{-1}+\lambda E_2(\lambda).
\end{align*}
The term \eqref{eq:S25} has been dealt with in Corollary~\ref{cor:P0}.  Now, we consider (\ref{eq:S24}). Note that
when we plug $R_V$ into (\ref{kxy}), then the term corresponding to the first summand in (\ref{eq:S24}) is
(with the notation $S=S_2b(0)^{-1}S_2$, $a_1=|y-y_1|$ and $a_2=|x-x_1|+|y-y_1|$)
\begin{align*}
& \frac{-1}{\pi i} \int_{-\infty}^{\infty} \int_{\R^6} e^{it\lambda^2}\chi_{\lambda_0}(\lambda)
\frac{e^{i\lambda|x-x_1|}-1}{\lambda 4 \pi |x-x_1|}
\frac{e^{i\lambda|y-y_1|}}{4\pi |y-y_1|}v(x_1)S(x_1,y_1)v(y_1) d x_1 dy_1d\lambda\\
&= \frac{-1}{16\pi^3} \int_{-\infty}^{\infty}\int_{\R^6} e^{it\lambda^2}\chi_{\lambda_0}(\lambda)
\frac{\sin(\lambda a_2)-\sin(\lambda a_1)} {\lambda  |x-x_1||y-y_1|}
 v(x_1)S(x_1,y_1)v(y_1) d x_1 dy_1 d\lambda\\
 &= \frac{-1}{16\pi^3} \int_{-\infty}^{\infty}\int_{\R^6} e^{it\lambda^2}\chi_{\lambda_0}(\lambda) \int_{a_1}^{a_2} \cos(\lambda b) db
\frac{v(x_1)S(x_1,y_1)v(y_1)}{|x-x_1||y-y_1|}
 d x_1 dy_1d\lambda \\
 &=:t^{-1/2} F_{1,t}(x,y).
\end{align*}
Arguing as in (\ref{fourier}), we obtain
\begin{align*}
|F_{1,t}(x,y)|&= c \left|\int_{-\infty}^{\infty}\int_{\R^6} e^{iu^2/4t}  \int_{a_1}^{a_2}
[\widehat{\chi_{\lambda_0}}(u+b)+ \widehat{\chi_{\lambda_0}}(u-b)] db \frac{v(x_1)S(x_1,y_1)v(y_1)}{|x-x_1||y-y_1|}
 d x_1 dy_1du \right| \\
& \lesssim \int_{\R^6} \int_{a_1}^{a_2} \int_{-\infty}^{\infty} | \widehat{\chi_{\lambda_0}}(u+b)+ \widehat{\chi_{\lambda_0}}(u-b)|
\frac{|v(x_1)| |S(x_1,y_1)| |v(y_1)| }{ |x-x_1||y-y_1|} du ,\ db \,d x_1 dy_1\\
& \lesssim \| \widehat{\chi_{\lambda_0}}\|_1  \int_{\R^6} \frac{|v(x_1)| |S(x_1,y_1)| |v(y_1)|}{|y-y_1|}
 d x_1 dy_1  \\
& \lesssim 1
\end{align*}
This inequality holds independently of $t, x$ and $y$. Therefore,
$$
\sup_t \| F_{1,t}\|_{L^1\rightarrow L^\infty} \lesssim 1 \text{\ \
and\ \ }\lim_{t\to\infty} F_{1,t}(x,y)=c\int_{\R^6}
\frac{v(x_1)S(x_1,y_1)v(y_1)}{|y-y_1|}\, dx_1dy_1.
$$
The second summand in (\ref{eq:S24}) can be treated similarly.

\nopar
Now, we consider \eqref{eq:S23}; it can be written as
\[ \eqref{eq:S23}  = \lambda^{-1} R_0(\lambda^2)vE_3(\lambda)v
R_0(\lambda^2)
\]
with
\[ \lambda E_3(\lambda) := -\Gamma_1(\lambda) S_1 \Gamma_2(\lambda)
S_2 b(\lambda)^{-1} S_2 \Gamma_2(\lambda) S_1 \Gamma_1(\lambda) + S_2b(0)^{-1}S_2.
\]
Clearly, the terms resulting from $E_3$ resemble $K_3$ from~\eqref{eq:K3def}. However, we do not have an extra
$\lambda$ at our disposal, which implies that instead of \eqref{eq:K3est} we will only obtain a $t^{-\half}$ power.
The details are as follows:  if we plug $R_V$ into (\ref{kxy}), then the term corresponding to  (\ref{eq:S23}) is
(up to constants)
\[
\int_{-\infty}^{\infty} \int_{\R^6} e^{it\lambda^2}\chi_{\lambda_0}(\lambda)
\frac{e^{i\lambda|x-x_1|}}{|x-x_1|}
\frac{e^{i\lambda|y-y_1|}}{|y-y_1|}v(x_1)E_3(\lambda)(x_1,y_1)v(y_1)\, dx_1 dy_1d\lambda.
\]
By the arguments that lead from \eqref{eq:K3begin} to \eqref{eq:K3end}, we conclude that the
absolute value of this expression does not exceed
\[ |t|^{-\half} \int_{-\infty}^\infty
\left\|\,\left|\widehat{\chi_{\lambda_0}E_3 }(\xi) \right|\, \right\|_{L^2\rightarrow L^2}d\xi.
\]
uniformly in $x,y\in\R^3$.
To bound this integral, we use Lemma~\ref{L:convolve}. Write
\begin{align*}
 E_3(\lambda) &= -\lambda^{-1}(\Gamma_1(\lambda)-\Gamma_1(0)) S_1 \Gamma_2(\lambda)
S_2 b(\lambda)^{-1} S_2 \Gamma_2(\lambda) S_1 \Gamma_1(\lambda) \\
& - S_1 \lambda^{-1}(\Gamma_2(\lambda)-\Gamma_2(0))
S_2 b(\lambda)^{-1} S_2 \Gamma_2(\lambda) S_1 \Gamma_1(\lambda) \\
& -  S_2 \lambda^{-1}(b(\lambda)^{-1}-b(0)^{-1}) S_2 \Gamma_2(\lambda) S_1 \Gamma_1(\lambda) \\
& -  S_2 b(0)^{-1} S_2 \lambda^{-1}(\Gamma_2(\lambda)-\Gamma_2(0)) S_1 \Gamma_1(\lambda) \\
& -  S_2 b(0)^{-1} S_2  \lambda^{-1}(\Gamma_1(\lambda)-\Gamma_1(0)).
\end{align*}
Consequently, we need to prove the bound of Lemma~\ref{L:hammal} for the following basic building blocks (we dropped the 
subscript $\lambda_0$):
\begin{align*}
 F_1(\lambda) &= \chi(\lambda)\Gamma_1(\lambda)=\chi(\lambda)(A(\lambda)+S_1)^{-1} \\
 F_2(\lambda) &= \chi(\lambda)\lambda^{-1}(\Gamma_1(\lambda)-\Gamma_1(0))=\chi(\lambda)\lambda^{-1}((A(\lambda)+S_1)^{-1}-(A_0+S_1)^{-1}) \\
 F_3(\lambda) &= \chi(\lambda)S_1\Gamma_2(\lambda)S_1=\chi(\lambda)S_1(m(\lambda)+S_2)^{-1}S_1 \\
 F_4(\lambda) &= \chi(\lambda)\lambda^{-1}S_1(\Gamma_2(\lambda)-\Gamma_2(0))S_1
 =\chi(\lambda)\lambda^{-1}S_1((m(\lambda)+S_2)^{-1}-(m(0)+S_2)^{-1})S_1
\end{align*}
as well as
\begin{align*}
F_5(\lambda) &= \chi(\lambda) S_2 b(\lambda)^{-1} S_2 = \chi(\lambda) S_2 (b(0)+\lambda b_1(\lambda))^{-1} S_2 \\
F_6(\lambda) &= \chi(\lambda) S_2 \lambda^{-1}(b(\lambda)^{-1}-b(0)^{-1}) S_2.
\end{align*}
The functions $F_j$ with $1\le j\le 4$ were already discussed in
Lemma~\ref{L:hammal}. The only difference here is the appearance
of $S_2$ in $F_3$ and $F_4$ (for the function $E_2$
see~\eqref{eq:m0inv}). But this does not effect the bounds from
Lemma~\ref{L:hammal}, which implies that we only need to prove the
following claims concerning the new terms $F_5$ and $F_6$: \beeq
\label{eq:F56} \max_{j=5,6}\int_{-\infty}^\infty \left\|\,
\left|\widehat{F_j }(\xi)\right|\, \right\|_{ 2\rightarrow 2}d\xi
< \infty. \eneq Recall that, see \eqref{eq:b1def},
\begin{align}
 b(0) &= -S_2 v G_2 v S_2 \nn \\
 b(\lambda) &= b(0) + \lambda b_1(\lambda) = b(0)(1+\lambda b(0)^{-1}b_1(\lambda)) \nn \\
 b_1(\lambda) &=  \f{S_2 [m_1(\lambda)-m_1(0)]S_2}{\lambda} +
\f{1}{\lambda} \sum_{k=1}^\infty (-1)^k\lambda^k S_2\left(m_1(\lambda) (m(0)+S_2)^{-1}\right)^{k+1}S_2 \label{eq:m1ser}\\
 b(\lambda)^{-1} &= \sum_{j=0}^\infty (-1)^j \lambda^j (b(0)^{-1}b_1(\lambda))^j b(0)^{-1}. \label{eq:b1ser}
\end{align}
Applying Lemma~\ref{L:convolve}  to the Neuman series in~\eqref{eq:b1ser} shows that in order to obtain~\eqref{eq:F56}, 
we need to prove that
\[ \int_{-\infty}^\infty \left\|\, \left|\widehat{\chi_1 b_1 }(\xi)\right|\,\right\|_{ 2\rightarrow 2}d\xi < \infty. \]
Another application of Lemma~\ref{L:convolve}, this time to the
Neuman series~\eqref{eq:m1ser}, reduces matters to proving
\[ \int_{-\infty}^\infty \left\|\, \left|\widehat{\chi_2 m_1 }(\xi)\right|\,\right\|_{ 2\rightarrow 2}d\xi < \infty, \]
which was already done in \eqref{eq:m1FT}. In both these cases,
the cut-off functions $\chi_1, \chi_2$ need to be taken with
sufficiently small supports. This leaves the term
\[ \f{S_2 [m_1(\lambda)-m_1(0)]S_2}{\lambda} \]
from~\eqref{eq:m1ser} to be considered.  In view of
\eqref{eq:m1def} and \eqref{eq:A2def},
\begin{align*}
& S_2\frac{m_1 (\lambda)-m_1(0)}{\lambda}S_2 \\
& = S_2
\f{A_2 (\lambda) - A_2 (0)}{\lambda}S_2
+\sum_{k=1}^{\infty} (-1)^k \lambda^{k-1 } S_2
\left(A_1(\lambda) (A_0 + S_1)^{-1}\right)^{k+1}S_2.
\end{align*}
By \eqref{chia1}, and Lemma~\ref{L:convolve}, the Neuman series
makes a summable contribution to~\eqref{eq:F56}. On the other
hand, the contribution of
\[ S_2\f{A_2 (\lambda) - A_2 (0)}{\lambda}S_2\]
to \eqref{eq:F56} is controlled by the bound~\eqref{a2prime}, and we are done.
\end{proof}


\end{document}